\documentclass[a4paper,10pt]{article}
\usepackage{amsmath}
\usepackage{amssymb,amsfonts,amsthm}
\usepackage[applemac]{inputenc}
\usepackage[english]{babel}
\usepackage[pdftex]{graphicx}
\usepackage{caption}
\usepackage{subcaption}
\usepackage{anysize}
\usepackage{enumerate}
\usepackage{tikz}
\usepackage[all]{xy}
\usepackage{mathrsfs}
\usepackage{verbatim}
\usepackage{cite}
\definecolor{green}{rgb}{0,0.8,0.5}

\newcommand{\abs}[1]{\left\vert#1\right\vert}

\graphicspath{{/pictures}}

\usepackage[pdftex]{hyperref}
\hypersetup{colorlinks=true,linkcolor=blue,citecolor=red}
\usepackage[width=0.8\textwidth]{caption}

\renewenvironment{abstract}{\small\quotation\noindent
 {\bfseries \abstractname .}}{\endquotation \par}

\newenvironment{trivlistparm}[2]{\trivlist
\item[{#1 #2}]}{\endtrivlist}

\newenvironment{demo}{\trivlistparm{\bfseries}{Proof.}}{\Qed\endtrivlistparm}

\catcode`\@=11

\def\resetthefootnote{\renewcommand{\thefootnote}{\@arabic\c@footnote} }
\def\@principiremex#1{\trivlist
 \item[\hskip \labelsep{\bfseries #1\ \thethm.}]\ignorespaces}
\def\opar@principiremex#1[#2]{\trivlist
 \item[\hskip \labelsep{\bfseries #1\ \thethm\ (#2).}]\ignorespaces}

\newcommand{\newTHEOremrom}[2]{\newenvironment{#1}{\refstepcounter{thm}\@ifnextchar[{\opar@principiremex{#2}}
{\@principiremex{#2}}}{\qedB\endtrivlist}} \catcode`\@=12
\DeclareMathSymbol{\square}{\mathord}{AMSa}{"03}
\newcommand{\qedB}{\nopagebreak\hspace*{\fill}$\square$\par}
\newcommand{\Qed}{\nopagebreak\hspace*{\fill}{\vrule width6pt height6pt depth0pt}\par}

\newTHEOremrom{defi} {Definition}
\newTHEOremrom {rem} {Remark}
\newTHEOremrom {ex} {Example}

\renewcommand{\epsilon}{\varepsilon}
\renewcommand{\leq}{\leqslant}
\newcommand{\R}{\mathbb{R}}

\newcommand{\D}{D_R}

\newcommand{\PA}{\mathscr{P}}

\newtheorem{thm}{Theorem}[section]

\newtheorem{cor}[thm]{Corollary}
\newtheorem{lema}[thm]{Lemma}

\title{\textbf{Bifurcation of relative equilibria generated by a circular vortex path in a circular domain}
\footnotetext{All the authors are partially supported by the MEC/FEDER grant MTM2014-52232-P.}
\footnotetext{{\it Email addresses:} \texttt{rojas@ugr.es}  (D.~Rojas, corresponding author), \texttt{ptorres@ugr.es} (P.J.~Torres).}
}

\author{D. Rojas and P.J. Torres\\[10pt]
{\small \textsl{Departamento de Matem\'atica Aplicada,}}\\
\vspace{-2pt}
{\small \textsl{Universidad de Granada, 18071 Granada, Spain}}}
\date{}

\begin{document}

\maketitle

\begin{abstract} We study the passive particle transport generated by a circular vortex path in a 2D ideal flow confined in a circular domain. Taking the strength and angular velocity of the vortex path as main parameters, the bifurcation scheme of relative equilibria is identified. For a perturbed path, an infinite number of orbits around the centers are persistent, giving rise to periodic solutions with zero winding number.
\end{abstract}

\section{Introduction}

The passive particle transport in a 2D incompressible flow with prescribed vorticity is a research topic of the highest relevance in Fluid Dynamics \cite{Aref2011,Saf,Ot}. In the Lagrangian formulation, the advection of single particles is ruled by a Hamiltonian system where the stream function plays the role of the Hamiltonian. In this paper, we consider the dynamics induced in an ideal flow confined in a circular domain of radius $R$ under the action of a prescribed $T$-periodic vortex path. This problem is classical in the literature (see for instance \cite{Aref84,BT13,CM,FrZ93}) and can be seen as a 2D idealization of the mixing of a fluid in a cylindrical tank.

Let $B_R\subset \R^2$ be the open ball of center $(0,0)$ and radius $R$, and consider a $T$-periodic vortex path given by $z:\R\to B_R$. Then, the stream function of the fluid confined in $B_R$ and under the action of the vortex is given by 
$$
\Psi(t,\zeta)=\frac{\Gamma}{2\pi}\left(\ln\abs{\zeta-z(t)}-\ln \abs{\zeta-\frac{R^2}{|z(t)|^2}z(t)}\right).
$$
Here, $\Gamma$ is the strength or charge of the vortex, and its sign gives the sense of rotation. In this function, the first term accounts for the vortex action, whereas the second term models the influence of the solid circular boundary. It is useful to see $\zeta$  as a complex variable, then the corresponding Hamiltonian system is 
\begin{equation}\label{vortex_eq}
\dot{\zeta^*} = \frac{\Gamma}{2\pi i}
\left(\frac{1}{\zeta-z(t)} - \frac{1}{\zeta - \tfrac{R^2}{\vert
z(t) \vert^2}z(t)}\right),
\end{equation}
where the asterisk means the complex conjugate. 

In the related literature, $z(t)$ is called the {\it stirring protocol}. When it is constant, then $\Psi$ is a conserved quantity and all the particles rotate around the vortex in circular trajectories. When $z(t)$ is time-dependent, then the Hamiltonian ceases to be a conserved quantity and the analysis is more delicate. In \cite{BT13}, it is proved that any smooth stirring protocol $z(t)$ induces an infinite number of periodic trajectories rotating around the vortex (non-zero winding number). Hence, it is a natural question to try to identify the stirring protocols that generate periodic trajectories with zero winding number, that is, particles moving periodically but not rotating around the vortex. This question was posed explicitly as an open problem in \cite[Section 7.3]{Torres}. Our intention is to advance on the comprehension of this difficult problem by analyzing the family of circular protocols $z(t)=r_0\exp(i\theta_0 t)$. In this case, the change to a corotating frame $\zeta(t)=\eta(t)\exp(i\theta_0 t)$ transforms \eqref{vortex_eq} into the autonomous system
\begin{equation}\label{vortex_eq2}
\dot{\eta}^* = i\theta_0\eta^*+\frac{\Gamma}{2\pi i}
\left(\frac{1}{\eta-r_0} - \frac{1}{\eta -
\tfrac{R^2}{r_0}}\right).
\end{equation}
The hamiltonian structure is preserved, so the streamlines are just the level
curves of corresponding Hamiltonian, which is indeed a conserved quantity. This fact enables a complete bifurcation analysis of the relative equilibria, performed in Section \ref{sec2}, which correspond to periodic solutions of zero winding number for the original system. Moreover, in Section \ref{sec3} it is shown that the bifurcation scheme identified in Section \ref{sec2} is robust under small perturbations of the vortex path. Finally, the last section exposes some conclusions of the presented study.

\begin{figure}
\begin{subfigure}{.24\textwidth}
\centering
\includegraphics[scale=.4]{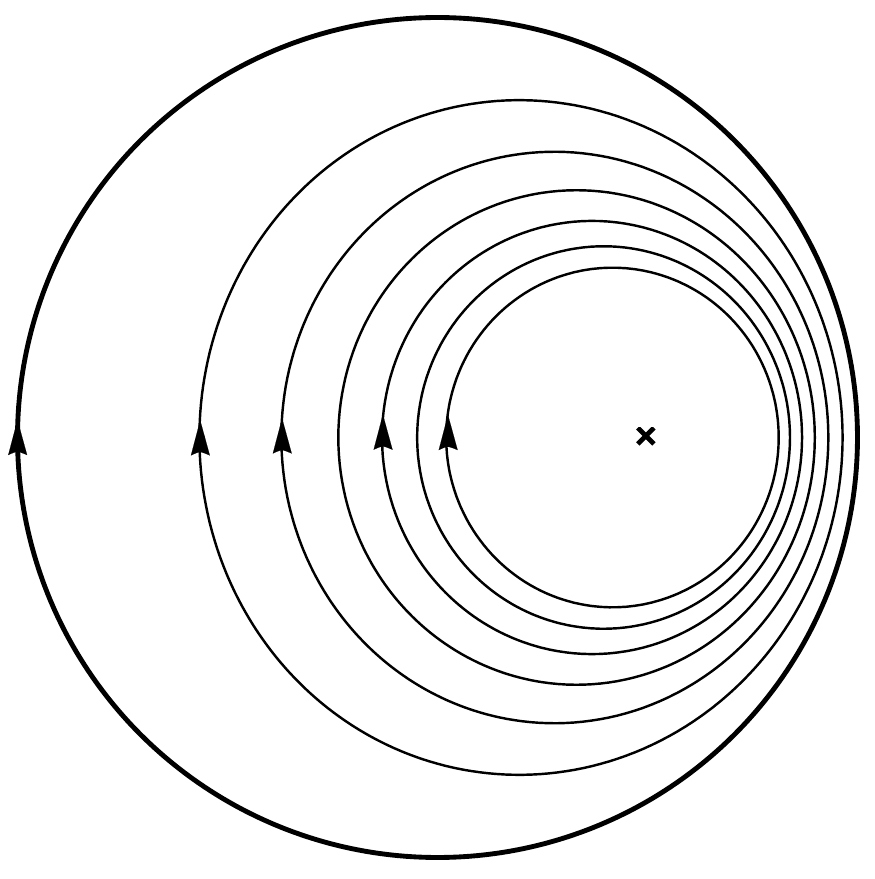}
\caption{\label{fig:retratos-a}}
\end{subfigure}
\begin{subfigure}{.24\textwidth}
\centering
\includegraphics[scale=.4]{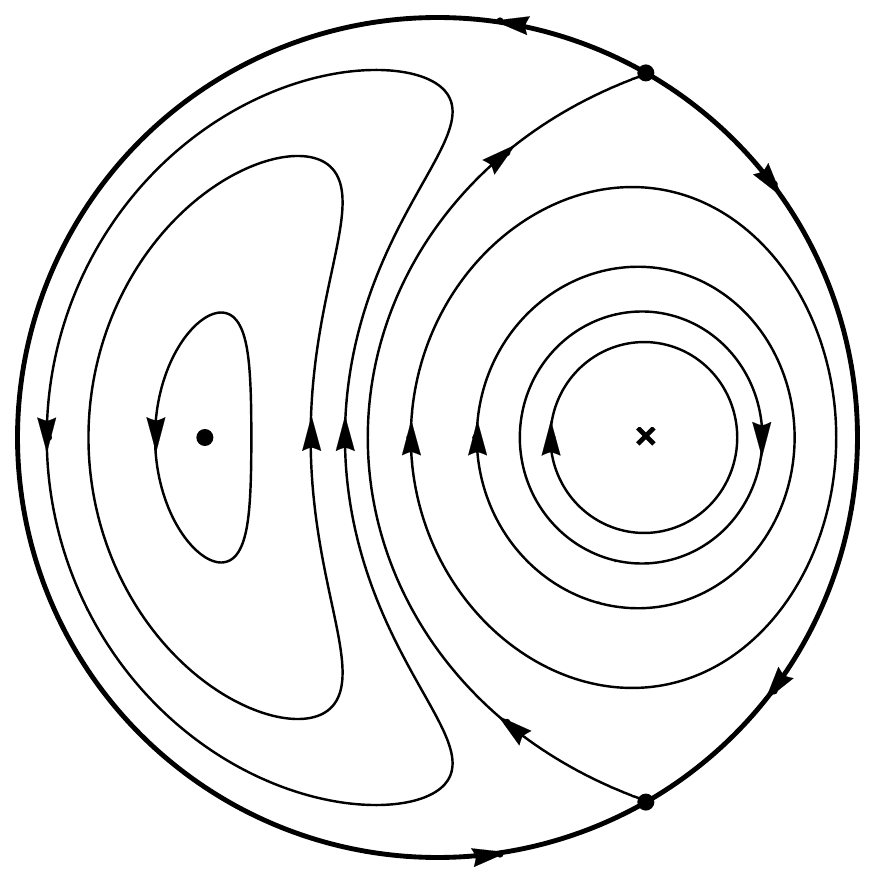}
\caption{\label{fig:retratos-b}}
\end{subfigure}
\begin{subfigure}{.24\textwidth}
\centering
\includegraphics[scale=.4]{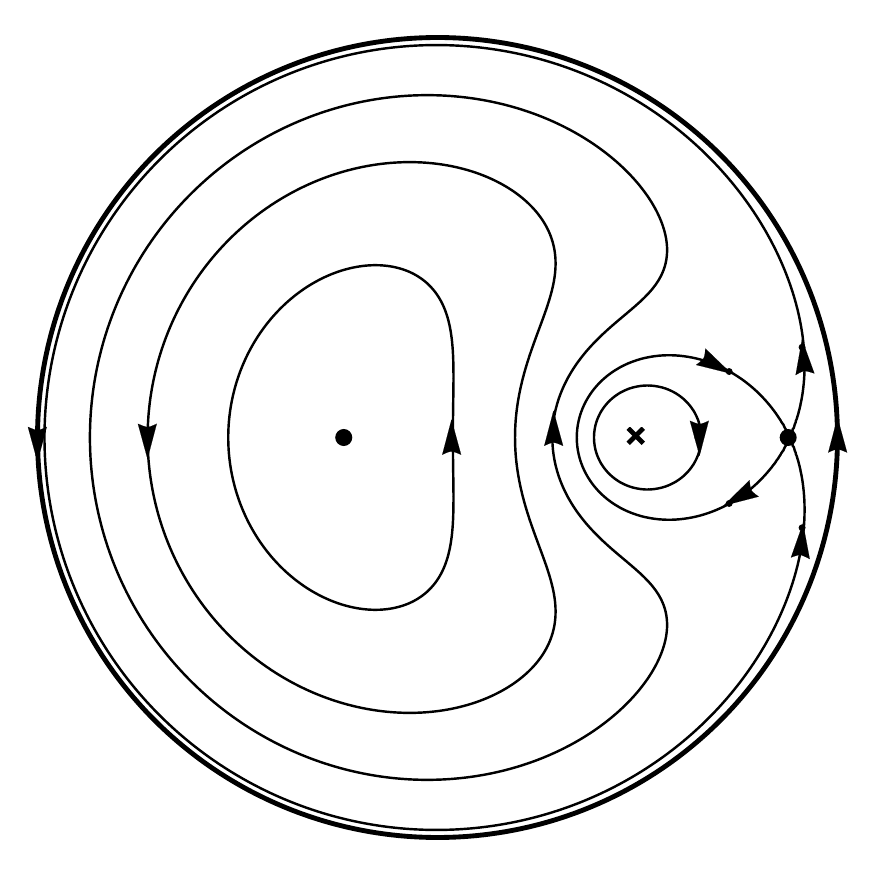}
\caption{\label{fig:retratos-c}}
\end{subfigure}
\begin{subfigure}{.24\textwidth}
\centering
\includegraphics[scale=.4]{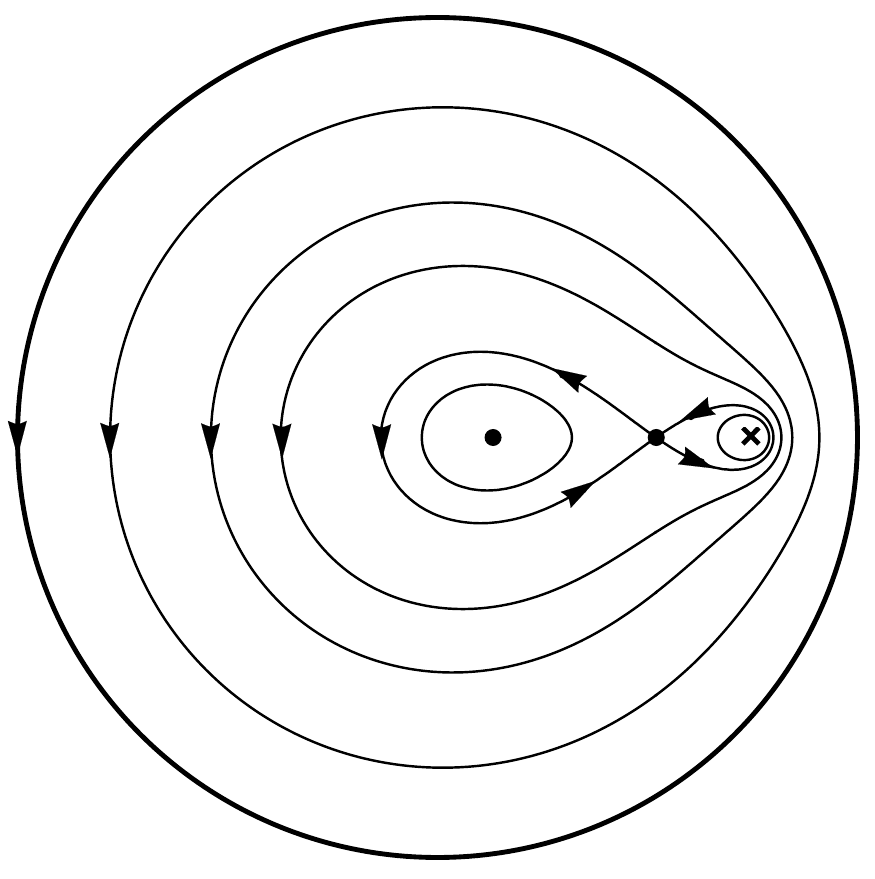}
\caption{\label{fig:retratos-d}}
\end{subfigure}
\begin{subfigure}{.24\textwidth}
\centering
\includegraphics[scale=.4]{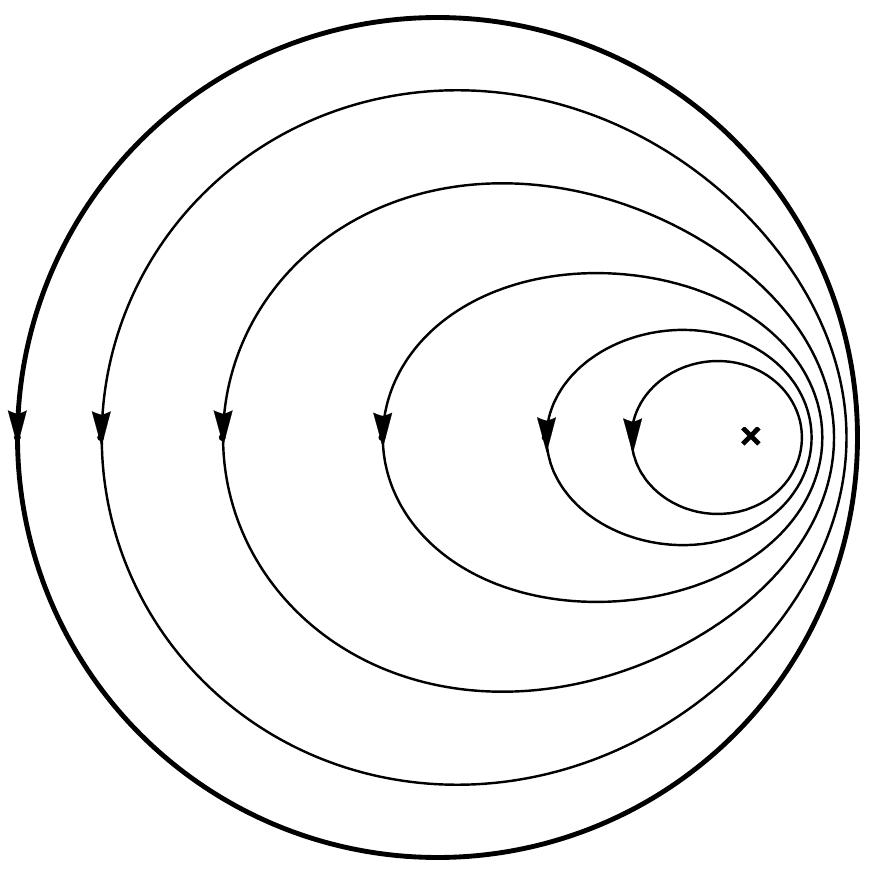}
\caption{\label{fig:retratos-e}}
\end{subfigure}
\begin{subfigure}{.24\textwidth}
\centering
\includegraphics[scale=.4]{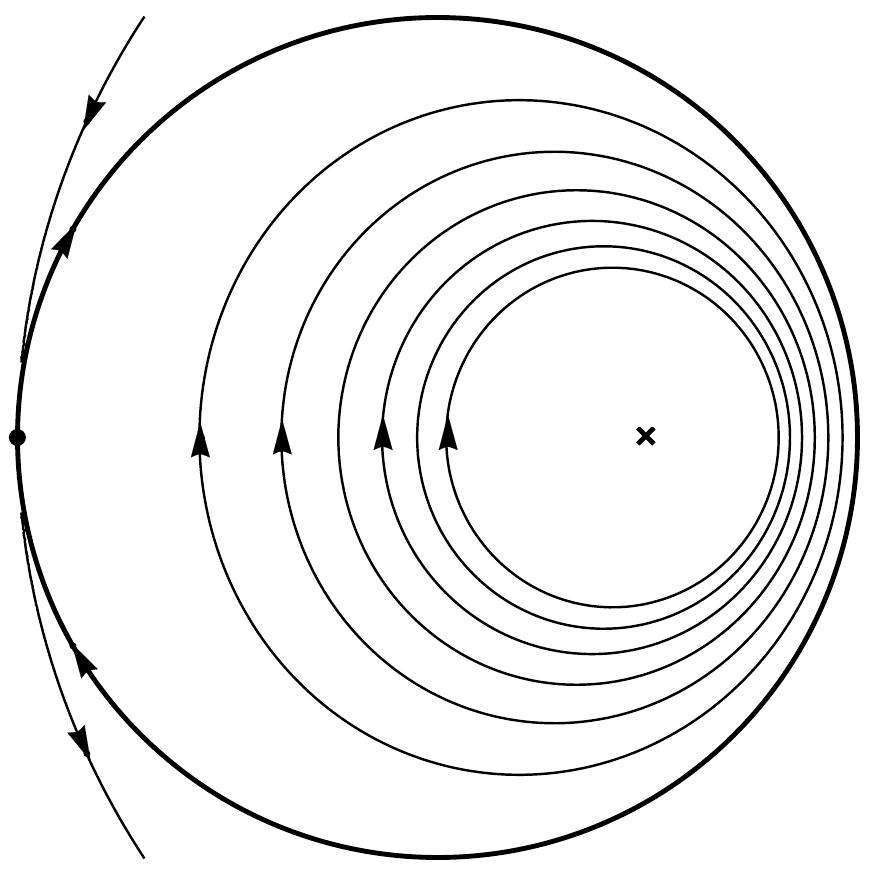}
\caption{\label{fig:retratos-f}}
\end{subfigure}
\begin{subfigure}{.24\textwidth}
\centering
\includegraphics[scale=.4]{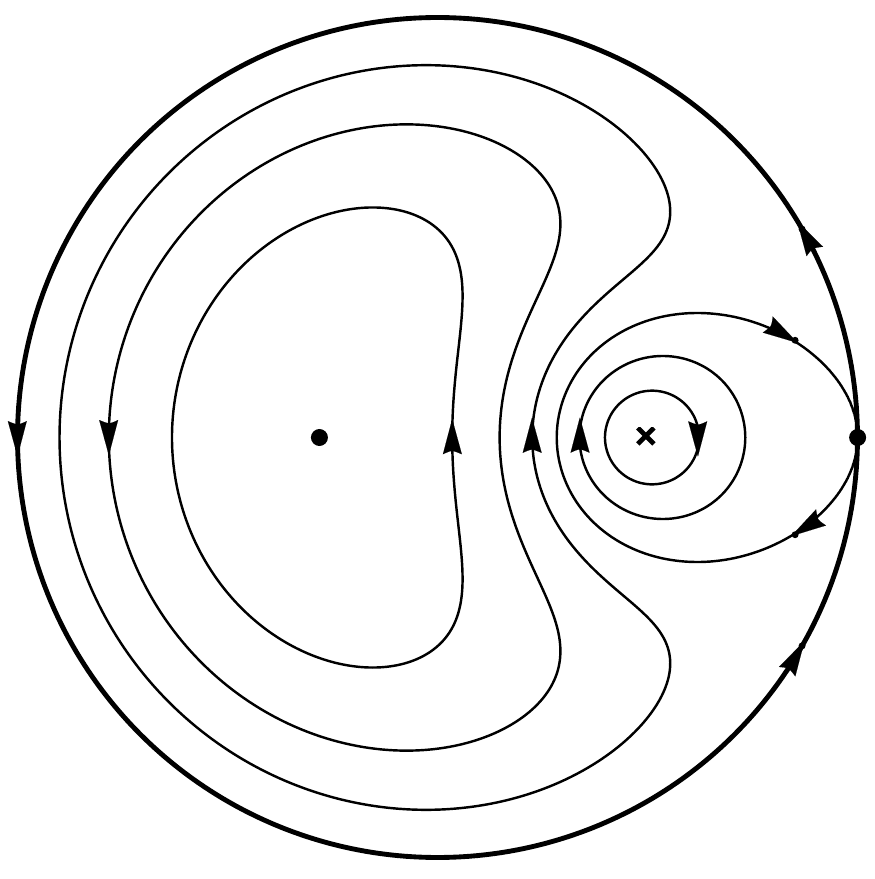}
\caption{\label{fig:retratos-g}}
\end{subfigure}
\begin{subfigure}{.24\textwidth}
\centering
\includegraphics[scale=.4]{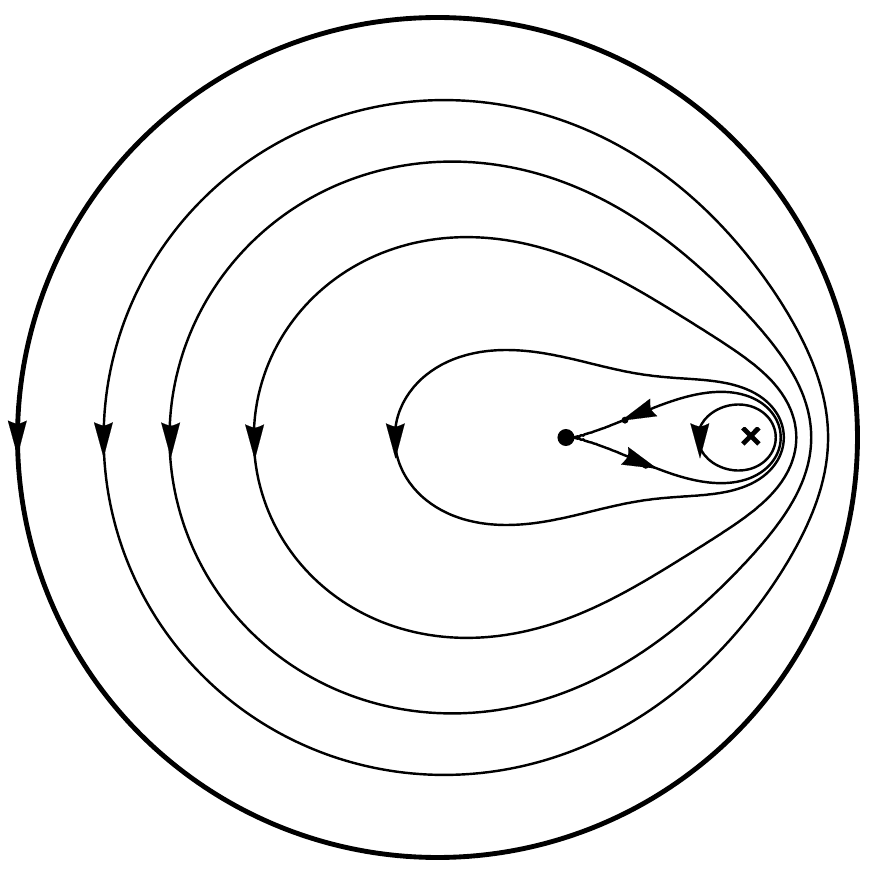}
\caption{\label{fig:retratos-h}}
\end{subfigure}
\caption{\label{fig:retratos}Phase portrait of system \eqref{eq:system} depending on the parameters according to Theorem~\ref{thm:retratos}.}
\end{figure}

\section{Phase portrait and bifurcations}\label{sec2}  

This section is devoted to the bifurcation analysis of the phase portrait of system \eqref{vortex_eq2}. Working on cartesian coordinates, the streamlines are level curves of the Hamiltonian function 
\[
\Psi(x,y)=-\frac{\theta_0}{2}(x^2+y^2)+\frac{\Gamma}{2\pi}\ln\sqrt{\frac{(x-r_0)^2+y^2}{(x-\frac{R^2}{r_0})^2+y^2}}.
\]
Here $(R,\Gamma,\theta_0)\in(0,+\infty)\times (\R\setminus\{0\})^2$ and $r_0\in(0,R)$. From now on, for the sake of further simplicity, we denote 
\[
a(x)\!:=a(x,r_0)=x-r_0,\  b(x)\!:=b(x,R,r_0)=x-\frac{R^2}{r_0} \text{ and }c\!:=\frac{\Gamma}{2\pi\theta_0}.
\] 
Thus, system \eqref{vortex_eq2} can be written in the $(x,y)$-variables as
\begin{equation}\label{eq:system}
\begin{cases}
\dot x = \dfrac{\partial \Psi}{\partial y}
=-\theta_0y+c\theta_0y\left( \dfrac{1}{a(x)^2+y^2} - \dfrac{1}{b(x)^2+y^2} \right),\\
\dot y = -\dfrac{\partial \Psi}{\partial x}
=\theta_0x-c\theta_0\left( \dfrac{a(x)}{a(x)^2+y^2} - \dfrac{b(x)}{b(x)^2+y^2} \right).\\
\end{cases}
\end{equation}

\begin{figure}
\centering
\includegraphics[scale=1]{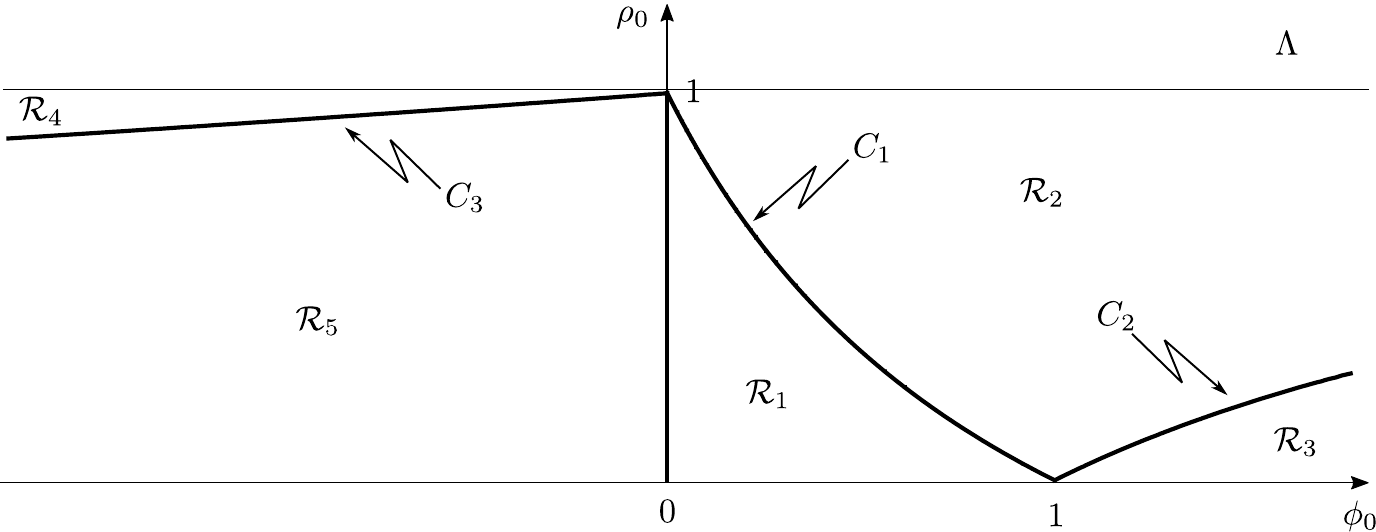}
\caption{\label{fig:diagrama-bif}Bifurcation diagram of the phase-portrait of system \eqref{eq:system} on $\D$. The blond curve correspond to bifurcation parameters $\mathcal B$, whereas the remaining ones correspond to regular parameters. In Theorem~\ref{thm:retratos} the phase portrait at each region is given.}
\end{figure}

Let $\D\subset\R^2$ be the closed ball of center $(0,0)$ and radius $R$. It is an immediate calculation to show that $\D$ is invariant by the flow of system \eqref{eq:system}. Next result deals with the phase portrait of the system on $\D$. It will be shown that the position $r_0$ and angular velocity $\theta_0$ of the vortex path are the main parameters on the system, whereas the remaining ones can be normalized. To this end, and for the sake of simplicity on the statement, we set $\rho_0\!:=\frac{r_0}{R}$ and $\phi_0\!:=\frac{R^2}{c}=\frac{2\pi R^2\theta_0}{\Gamma}$. Thus, the parameter space of system~\eqref{eq:system} turns $\Lambda\!:=\{(\rho_0,\phi_0)\in\R^2: 0<\rho_0<1 \text{ and }\phi_0\neq 0\}$. Moreover, let us define
\[
f(\rho_0,\phi_0)\!:=27\rho_0^2(\rho_0^2-1)+\phi_0\left(2-3\rho_0^2-3\rho_0^4+2\rho_0^6-2(1-\rho_0^2+\rho_0^4)^{\frac{3}{2}}\right),
\]
and
\[
\mathcal B\!:=\left\{(\rho_0,\phi_0)\in\Lambda : \phi_0f(\rho_0,\phi_0)\!\left(\rho_0-\frac{1-\phi_0}{1+\phi_0}\right)\!\!\left(\rho_0-\frac{\phi_0-1}{1+\phi_0}\right)=0\right\}.
\]
The curve $\mathcal B$ is the union of three curves, namely $C_i$, $i=1,\dots,3$ and splits the parameter space $\Lambda$ into five connected components, $\mathcal R_i$, $i=1,\dots,5$, according with Figure~\ref{fig:diagrama-bif}.   

\begin{thm}\label{thm:retratos}
Let $(\rho_0,\phi_0)\in \Lambda$. The set $\Lambda\setminus\mathcal B$ corresponds to regular parameters of system~\eqref{eq:system}. On each connected component, the phase portrait is the following:
\begin{enumerate}[$(a)$]
\item If $(\rho_0,\phi_0)\in\mathcal R_1$ then the dynamics on $\D$ is a global vortex at $(r_0,0)$ (see Figure~\ref{fig:retratos-a}).

\item If $(\rho_0,\phi_0)\in\mathcal R_2$ then the system has a vortex at $(r_0,0)$, a center at $(x_c^*,0)$ with $x_c^*\in(-R,0)$ and two hyperbolic saddles $(x_{s}^*,\pm y_{s}^*)$ at $\partial\D$ with a saddle connection inside $\D$ (see Figure~\ref{fig:retratos-b}).

\item If $(\rho_0,\phi_0)\in\mathcal R_3$ then the system has a vortex at $(r_0,0)$, a center at $(x_c^*,0)$ with $x_c^*\in(-R,0)$ and a hyperbolic saddle at $(x_s^*,0)$ with $x_s^*\in(r_0,R)$ (see Figure~\ref{fig:retratos-c}).

\item If $(\rho_0,\phi_0)\in\mathcal R_4$ then the system has a vortex at $(r_0,0)$, a center $(x_c^*,0)$ and a hyperbolic saddle $(x_s^*,0)$ satisfying $0<x_c^*<x_s^*<r_0$ (see Figure~\ref{fig:retratos-d}).

\item If $(\rho_0,\phi_0)\in\mathcal R_5$ then the dynamics on $\D$ is a global vortex at $(r_0,0)$ (see Figure~\ref{fig:retratos-e}).
\end{enumerate}
Moreover, the set $\mathcal B$ corresponds to bifurcation parameters of system~\eqref{eq:system}. On each curve the phase portrait is the following:
\begin{enumerate}[$(a)$]\setcounter{enumi}{5}
\item If $(\rho_0,\phi_0)\in C_1$ then the system has a vortex at $(r_0,0)$ and a degenerated saddle at $(-R,0)$ (see Figure~\ref{fig:retratos-f}). 

\item If $(\rho_0,\phi_0)\in C_2$ then the system has a vortex at $(r_0,0)$, a center at $(x_c^*,0)$ with $x_c^*\in(-R,0)$ and a degenerated saddle at $(R,0)$ (see Figure~\ref{fig:retratos-g}).

\item If $(\rho_0,\phi_0)\in C_3$ then the system has a vortex at $(r_0,0)$ and a cusp at $(x_{p}^*,0)$ (see Figure~\ref{fig:retratos-h}), where 
\[
x_{p}^*\!:=x_{p}^*(R,r_0)=\frac{R^2+r_0^2-\sqrt{R^4-R^2r_0^2+r_0^4}}{3r_0}.
\]
\end{enumerate}
\end{thm}

\begin{demo}
For the sake of simplicity we first begin the proof by showing that the only critical points in $\D$ that do not lie on the line $\{y=0\}$ are the hyperbolic saddles $(x_{s}^*,\pm y_{s}^*)$ at $\partial\D$ of case $(b)$ on the statement. To this end, assuming $y\neq 0$, from equations in \eqref{eq:system} we have that $\dot x=0$ if and only if
\[
-1+c\left(\frac{1}{a(x)^2+y^2}-\frac{1}{b(x)^2+y^2}\right)=0.
\]
Since $a(x)^2<b(x)^2$ for all $x<R$, if $\phi_0<0$ then $c<0$ and so the left hand side of the previous equality is negative. Then assume $\phi_0>0$. In this case, $\dot x=0$ if and only if
\[
y^2=-\frac{1}{2}(a(x)^2+b(x)^2)+\frac{1}{2}\sqrt{(b(x)^2-a(x)^2)(4c+b(x)^2-a(x)^2)}.
\]
Substituting the previous equality on the expression of $\dot y$ in \eqref{eq:system} and equaling zero one gets the equation
\[
\frac{(2x-a(x)-b(x))(a(x)+b(x))+\sqrt{(b(x)^2-a(x)^2)(4c+b(x)^2-a(x)^2)}}{2(a(x)+b(x))}=0.
\]
Thus, using that $a(x)=x-r_0$ and $b(x)=x-\frac{R^2}{r_0}$, the previous equation has the unique solution 
\[
x^*_{s}=\frac{R^2+r_0^2}{2r_0}-\frac{c}{2r_0}\left(1-\frac{r_0^2}{R^2}\right)
\]
and so
\[
(y^*_{s})^2=\frac{1}{4}\left(\frac{2(c^2+R^4)}{R^2}-\frac{(c-R^2)^2}{r_0^2}-\frac{(c+R^2)^2r_0^2}{R^4}\right).
\]
It is a computation to show that $x_{s}^*\in(-R,R)$ if and only if $\rho_0>\max\{\frac{1-\phi_0}{1+\phi_0},\frac{\phi_0-1}{1+\phi_0}\}$ (and so if and only if $(\rho_0,\phi_0)\in\mathcal R_2$) and $(x_{s}^*)^2+(y_{s}^*)^2=R^2$. It is only remaining to prove that $(x_{s}^*,\pm y_{s}^*)$ are hyperbolic saddles. This can be done evaluating the previous expression of the points $(x_{s}^*,\pm y_{s}^*)$ on the Jacobian matrix of system \eqref{eq:system}. In the case of $(x_{s}^*,y_{s}^*)$ the determinant of the Jacobian matrix is
\[
\text{det}(DX(x_{s}^*,y_{s}^*))=\frac{(cR-R^3-(c+R^2)r_0)(cR-R^3+(c+R^2)r_0)\theta_0^2}{c^2(R^2-r_0^2)}
\]
which is negative if and only if $\rho_0>\max\{\frac{1-\phi_0}{1+\phi_0},\frac{\phi_0-1}{1+\phi_0}\}$. Then, $(x_{s}^*,y_{s}^*)$ is a hyperbolic saddle. The same argument is valid for $(x_{s}^*,-y_{s}^*)$. Moreover, since $\partial\D$ is an invariant curve of system~\eqref{eq:system} and $(x_{s}^*,y_{s}^*)\in\partial\D$, $\partial\D$ is the stable manifold of one saddle (and unstable of the other). The correspondent unstable (stable) manifold cuts transversally the disk of radius $R$ due to the hyperbolicity of the saddles and so the connection between the saddles follows by Poincar\'e-Bendixon's theorem. 

The previous argument shows that out of case $(b)$ on the statement, all the critical points of system \eqref{eq:system} lie on $\{y=0\}$. Let us prove now the remaining cases of the result. Let us first consider $\phi_0>0$. This correspond to statements $(a)-(c)$ and $(f)-(g)$.
We can assume with no loss of generality that $\theta_0>0$ and $\Gamma>0$. The case with $\theta_0<0$ and $\Gamma<0$ follows by reversion of time. Notice that the hypothesis $\phi_0>0$ implies $c>0$. System \eqref{eq:system} has critical point at $(x^*,0)$ inside the disk of radius $R$ if and only if the function
\[
F(x)\!:=\theta_0\left(x-c\left(\frac{1}{a(x)}-\frac{1}{b(x)}\right)\right)
\]
satisfies $F(x^*)=0$ for some $x^*\in(-R,R)$. Multiplying by $a(x)b(x)$ the previous condition turns into $F(x^*)a(x^*)b(x^*)=0$. We point out that, on account of the expressions of $a(x)$ and $b(x)$, the previous two conditions are equivalent if $x^*\notin\{r_0,\frac{R^2}{r_0}\}$ (those correspond to singularities on the Hamiltonian function and so no critical points). Thus, system \eqref{eq:system} has a critical point at $(x^*,0)$ in $\D$ if and only if
\begin{equation}\label{eq:condition}
x^*a(x^*)b(x^*)=c\left(r_0-\frac{R^2}{r_0}\right)=:\!\lambda=\lambda(r_0,R,c).
\end{equation}
Notice that, since $\phi_0>0$, $\lambda<0$. The cubic polynomial $P(x)\!:=xa(x)b(x)$ has zeros at $x=0$, $x=r_0$ and $x=\frac{R^2}{r_0}$. $P(x)$ is negative if $x\in(-\infty,0)\cup(r_0,R^2/r_0)$ and it is positive if $x\in(0,r_0)\cup(R^2/r_0,+\infty)$, and the local maximum and minimum are, respectively,
\[
x_M=\frac{R^2+r_0^2-\sqrt{R^4-R^2r_0^2+r_0^4}}{3r_0}, \ \ x_m=\frac{R^2+r_0^2+\sqrt{R^4-R^2r_0^2+r_0^4}}{3r_0}.
\]
On the other hand, since $\phi_0>0$ then $\lambda=\lambda(r_0,R,c)$ varies from zero to $-\infty$. Thus $P(x)-\lambda=0$ has always a solution $x^*\in(-\infty,0)$, it has a double solution $x^*=x_m$ if $P(x_m)=\lambda$ and two solutions in $(r_0,R^2/r_0)$ if $P(x_m)<\lambda$. Let us study when this solutions correspond to critical points in $\D$. We point out that $x_m>R$ so at most two zero of $P(x)-\lambda$ lie in $(-R,R)$. It is a computation to show that $P(R)-\lambda>0$ if $\rho_0>\frac{\phi_0-1}{1+\phi_0}$, $P(R)-\lambda=0$ if $\rho_0=\frac{\phi_0-1}{1+\phi_0}$ and $P(R)-\lambda<0$ if $\rho_0<\frac{\phi_0-1}{1+\phi_0}$. On the other hand, $P(-R)-\lambda>0$ if $\rho_0<\frac{1-\phi_0}{1+\phi_0}$, $P(-R)-\lambda=0$ if $\rho_0=\frac{1-\phi_0}{1+\phi_0}$ and $P(-R)-\lambda<0$ if $\rho_0>\frac{1-\phi_0}{1+\phi_0}$.
Thus, if $(\rho_0,\phi_0)\in\mathcal R_1$, no roots of $P(x)-\lambda$ are inside $[-R,R]$ and so the result in $(a)$ holds. If $(\rho_0,\phi_0)\in C_1$, the unique zero of $P(x)-\lambda$ in $[-R,R]$ is $x=-R$. This correspond to a critical point of system \eqref{eq:system} at $(-R,0)$. Moreover, it is a degenerated saddle since $\partial\D$ is an invariant curve of system \eqref{eq:system} so $(f)$ is proved. If $(\rho_0,\phi_0)\in\mathcal R_2$, $P(x)-\lambda$ has only one zero $x_C^*\in(-R,0)$. Then, on account of the previous discussion about the hyperbolic saddles $(x_s^*,\pm y_s*)$ on $\partial\D$, result in $(b)$ is proved.
If $(\rho_0,\phi_0)\in C_2$, $P(x)-\lambda$ has two zeros: $x=x_c^*\in(-R,0)$ and $x=R$. The critical point $(R,0)$ correspond to a degenerated saddle since $\partial\D$ is an invariant curve of system \eqref{eq:system}. Then $(g)$ holds.
Finally, if $(\rho_0,\phi_0)\in\mathcal R_3$, $P(x)-\lambda$ has two zeros: $x=x_c^*\in(-R,0)$ and $x=x_s^*\in(r_0,R)$. In order to end with the case $\phi_0>0$ it only remains to prove that $x_c^*$ and $x_s^*$ are a center and a hyperbolic saddle, respectively.

The Jacobian matrix associated to system \eqref{eq:system} with $y=0$ is given by
\begin{equation}\label{eq:dif}
DX(x,0)=\left(
\begin{matrix}
0 & \theta_0\left(-1+c\left( \dfrac{1}{a(x)^2}-\dfrac{1}{b(x)^2} \right)\right)\\
\theta_0\left(1+c\left( \dfrac{1}{a(x)^2}-\dfrac{1}{b(x)^2} \right)\right) & 0\\
\end{matrix}\right).
\end{equation}
Notice that $\theta_0\left(1+c\left( \dfrac{1}{a(x)^2}-\dfrac{1}{b(x)^2} \right)\right)>0$ for all $x\in(-R,R)$. On the other hand, setting $x=x^*$ a critical point of \eqref{eq:system}, we have
\begin{align*}
c\left(\frac{1}{a(x^*)^2}-\frac{1}{b(x^*)^2}\right)-1
&=x^*\left(\frac{1}{a(x^*)}+\frac{1}{b(x^*)}\right)-1\\
&=\frac{(x^*)^2-R^2}{(x^*-r_0)(x^*-\frac{R^2}{r_0})},
\end{align*}
where we used $F(x^*)=0$ on the first equality and the expressions of $a(x)$ and $b(x)$ on the second. Thus
\begin{equation}\label{eq:dif2}
DX(x^*,0)=\left(
\begin{matrix}
0 & \theta_0\frac{(x^*)^2-R^2}{(x^*-r_0)(x^*-\frac{R^2}{r_0})}\\
\theta_0\left(2+\frac{(x^*)^2-R^2}{(x^*-r_0)(x^*-\frac{R^2}{r_0})}\right) & 0\\
\end{matrix}\right).
\end{equation}
Consequently, taking $x^*=x_c^*\in(-R,0)$ we have $(x_c^*)^2-R^2<0$ and $(x_c^*-r_0)(x_c^*-\frac{R^2}{r_0})>0$. Therefore $\frac{(x_c^*)^2-R^2}{(x_c^*-r_0)(x_c^*-\frac{R^2}{r_0})}<0$ and so $\det(DX(x_c^*,0))>0$. This implies that $(x_c^*,0)$ is a center. On the other hand, taking $x^*=x_s^*\in(r_0,R)$, $\frac{(x_s^*)^2-R^2}{(x_s^*-r_0)(x_s^*-\frac{R^2}{r_0})}>0$ and so $\det(DX(x_s^*,0))<0$. This implies that $(x_s^*,0)$ is a hyperbolic saddle. This ends with the proof of statements $(a)$, $(b)$, $(c)$, $(f)$ and $(g)$.

Let us now consider the case $\phi_0<0$. This corresponds to statements $(d)$, $(e)$ and $(h)$. In this situation we can assume with no loss of generality that $\theta_0>0$ and $\Gamma<0$. The opposite case follows by reversion of time. Notice that the hypothesis $\phi_0<0$ implies $c<0$. Consequently, on account of the equation \eqref{eq:system} critical points can only belong to $\{(x,y)\in\R^2:y=0\}$. Similarly as before, system \eqref{eq:system} has a critical point at $(x^*,0)$ inside the disk of radius $R$ if and only if \eqref{eq:condition} is satisfied. Notice that, since $\phi_0<0$, in this case $\lambda=\lambda(r_0,R,c)$ varies from zero to $+\infty$. Thus, on account of $R<R^2/r0$, if $\lambda$ stays above of the maxima of $P(x)$ inside $(0,r_0)$ then $P(x)-\lambda=0$ has a unique zero which is larger than $R$. This happens when $f(\rho_0,\phi_0)<0$. If $f(\rho_0,\phi_0)=0$ then the maximum of $P(x)$ inside $(0,r_0)$ contact $\lambda$ and gives the cusp $(x_{p}^*,0)$ with $x_{p}^*=x_M$. Finally, if $f(\rho_0,\phi_0)>0$ then the maximum of $P(x)$ is greater than $\lambda$ and so $P(x)-\lambda$ has two real roots inside $(0,r_0)$: namely $x_c^*$ and $x_s^*$, satisfying $0<x_c^*<x_M<x_s^*<r_0$. It only remains to prove the stability of such critical points. This follows from the expression in \eqref{eq:dif} of the Jacobian matrix associated to system \eqref{eq:system} with $y=0$. We point out that, since $a(x)^2<b(x)^2$ for all $x\in(0,r_0)$ and $c<0$, we have 
\[
\theta_0\left(-1+c\left( \dfrac{1}{a(x)^2}-\dfrac{1}{b(x)^2} \right)\right)<0.
\]
On the other hand, setting $x=x^*$ a critical point of system \eqref{eq:system}, on account of $F(x^*)=0$ we have
\[
1+c\left(\frac{1}{a(x^*)^2}-\frac{1}{b(x^*)^2}\right)=2+\frac{(x^*)^2-R^2}{(x^*-r_0)(x^*-\frac{R^2}{r_0})}=\frac{3r_0(x^*)^2-2(R^2+r_0^2)x^*+R^2r_0}{(r_0-x^*)(R^2-r_0x^*)}.
\]
The previous expression is positive if $x^*\in(0,x_M)$ and it is negative if $x^*\in(x_M,r_0)$. This proves that $x_{c}^*$ is a center and $x_s^*$ is an hyperbolic saddle and ends with the proof of $(d)$, $(e)$ and $(h)$.
\end{demo}

\section{Periodic perturbations and local continuation of periodic orbits}\label{sec3}

Given an autonomous planar Hamiltonian system 
\begin{equation}\label{eq:aut-H}
\dot \eta = J\nabla \mathcal H(\eta),
\end{equation}
it is interesting to ask about the existence of periodic solutions of the non-autonomous planar Hamiltonian system
\begin{equation}\label{eq:nonaut-H}
\dot \eta = J\nabla H(t,\eta;\epsilon),
\end{equation}
which are small $T$-periodic perturbations of \eqref{eq:aut-H} meaning that $H(t,\eta;0)\equiv\mathcal H(\eta)$. A.~Fonda, M.~Sabatini and F. Zanolin prove in \cite{FSZ} that under the hypothesis of the existence of a non-isochronous period annulus for the autonomous Hamiltonian system and some regularity conditions on $H(t,\eta;\epsilon)$ such periodic orbits exist.

More precisely, consider $\mathcal H:\mathcal A\rightarrow\R$ twice continuously differentiable and $\mathcal A\subseteq\R^2$ a period annulus such that the inner and outer components of its boundary are Jordan curves. Assume that $\mathcal A$ is not isochronous, that is the period of the periodic orbits in $\mathcal A$ covers an interval $[\mathcal T_{\min},\mathcal T_{\max}]$, with $\mathcal T_{\min}<\mathcal T_{\max}$. Then consider $H:\R\times\mathcal A\times(0,\epsilon_0)\rightarrow\R$, whose gradient with respect to the second variable, denoted by $\nabla H(t,\eta;\epsilon)$, is continuous in $(t,\eta;\epsilon)$, locally Lipschitz continuous in $\eta$ and $T$-periodic in $t$ for some $T>0$. Under these assumptions, the authors in \cite{FSZ} prove the following result:

\begin{thm}[Fonda, Sabatini and Zanolin]\label{thm:FSZ}
Given two positive integers $m$ and $n$ satisfying
\begin{equation}\label{T-condition}
\mathcal T_{\min}<\frac{mT}{n}<\mathcal T_{\max},
\end{equation}
there is an $\bar\epsilon>0$ such that, if $\abs{\epsilon}\leq\bar\epsilon$, then system \eqref{eq:nonaut-H} has at least two $mT$-periodic solutions, whose orbits are contained in $\mathcal A$, which make exactly $n$ rotations around the origin in the period time $mT$.
\end{thm}

The authors also emphasize the following immediate consequence:

\begin{cor}\label{cor:FSZ}
For any positive integer $N$ there is a $\bar\epsilon_N>0$ such that, if $\abs{\epsilon}<\bar\epsilon_N$, then system \eqref{eq:nonaut-H} has at least $N$ periodic solutions, whose orbits are contained in $\mathcal A$.
\end{cor}

Our propose in this section is to illustrate this situation in the case when system \eqref{eq:system} has a non-degenerated center inside $\D$. This occurs for parameters $(\rho_0,\phi_0)\in\mathcal R_2\cup\mathcal R_3\cup\mathcal R_4$, corresponding to the phase portraits $(b)$, $(c)$ and $(d)$ in Figure~\ref{fig:retratos} and Theorem~\ref{thm:retratos}. Let us denote by $\PA$ the period annulus of the center. The inner boundary of $\PA$ is the center itself, namely $p$, whereas the outer boundary of $\PA$ is formed by saddle connections. In both cases the outer boundary have critical points so it is clear that the period function tends to infinity as the orbits approach the outer boundary. Particularly, the center is not isochronous. Next result states the period of the linearized center.

\begin{lema}\label{lema:periodo-lineal}
Let $(\rho_0,\phi_0)\in\mathcal R_2\cup\mathcal R_3\cup\mathcal R_4$ and let $p=(x,0)$ be the non-degenerated center of system \eqref{eq:system}. Then the period of the associated linearized system at $p$ is
\[
T_0(x)=\frac{2\pi}{\sqrt{\theta_0^2\nu(\frac{x}{R})(2+\nu(\frac{x}{R}))}}
\]
where $\nu(x)\!:=\dfrac{x^2-1}{(x-\rho_0)(x-\frac{1}{\rho_0})}$.
\end{lema}

\begin{demo}
From the expression of the Jacobian matrix of the system in \eqref{eq:dif2} we have that the eigenvalues associated to the center are
\[
\lambda_{\pm}=\pm\omega i=\pm\sqrt{\theta_0^2 \frac{x^2-R^2}{(x-r_0)(x-\frac{R^2}{r_0})}\left(2+\frac{(x^2-R^2)}{(x-r_0)(x-\frac{R^2}{r_0})}\right)}i
\]
where $\omega$ denotes the frequency of the linearized center. Thus, setting $\rho_0=\frac{r_0}{R}$, and using that the period of the linearized center is $T_0=\frac{2\pi}{\omega}$ the result holds.
\end{demo}

Let us now consider the periodically perturbed stirring protocol $z_{\epsilon}(t)=r_{\epsilon}(t)\exp(i\theta_0 t)$ on system \eqref{vortex_eq} where $r_{\epsilon}(t)$ is a smooth $T$-periodic perturbation of $r_0$. More concretely, $r_{\epsilon}(t)=r_0+\epsilon f(t)+g(t;\epsilon)$ with $f$ and $g(\cdot;\epsilon)$ $T$-periodic analytic functions and $g(t;\epsilon)$ tending to zero uniformly on $t\in\R$ as $\epsilon$ tends to zero. The same change to a corotating frame than the presented at the beginning of this paper transforms \eqref{vortex_eq} into a periodic Hamiltonian system with Hamiltonian function

\begin{equation}\label{hamiltonian-perturbed}
\Psi(t,x,y;\epsilon)=-\frac{\theta_0}{2}(x^2+y^2)+\frac{\Gamma}{2\pi}\ln\sqrt{\frac{(x-r_{\epsilon}(t))^2+y^2}{(x-\frac{R^2}{r_{\epsilon}(t)})^2+y^2}}.
\end{equation}

\begin{thm}
Let $(\rho_0,\phi_0)\in\mathcal R_2\cup\mathcal R_3\cup\mathcal R_4$.
For any positive integer $N$ there is a $\bar\epsilon_N>0$ such that, if $\abs{\epsilon}<\bar\epsilon_N$, then the Hamiltonian system $\dot u = J\nabla \Psi(t,u;\epsilon)$ has at least $N$ periodic solutions contained in $D_R$. Particularly, the flow induced by system \eqref{vortex_eq} with $T$-periodic protocol $z_{\epsilon}(t)$ has infinity many periodic trajectories with zero winding number.
\end{thm}

\begin{demo}
The spirit of this proof is to use Theorem~\ref{thm:FSZ} in a certain period annulus where the regularity hypothesis are satisfied. On the one hand, since the outer boundary of the whole period annulus of the center in system~\eqref{eq:system} is a saddle connection, the period of the periodic orbits tends to infinity as they approach the outer boundary. On the other hand, setting $p=(x^*,0)$ the center itself, Lemma~\ref{lema:periodo-lineal} states that the period tends to
\[
T_0(x^*)=\frac{2\pi}{\sqrt{\theta_0^2\nu(\frac{x^*}{R})(2+\nu(\frac{x^*}{R}))}}
\]
as the orbits tends to $p$. The analyticity of the period function ensures then that for any $M>0$ large enough there exists a period annulus, namely $\mathcal A_M$, such that the period of its orbits covers $[T_0(x^*),M]$.

Let us apply Theorem~\ref{thm:FSZ} in $\mathcal A_M$. To this end, it is enough to show that $\nabla \Psi(t,\eta;\epsilon)$ is continuous in $(t,\eta;\epsilon)\in\R\times\mathcal A_M\times(0,\epsilon_0)$ for some $\epsilon_0>0$, locally Lipschitz continuous in $\eta=(x,y)\in\mathcal A_M$ and $T$-periodic in $t$. From the expression in \eqref{hamiltonian-perturbed},
\begin{equation}\label{grad-ham-per}
\nabla \Psi(t,\eta;\epsilon)=\left(
\begin{array}{c}
-y\theta_0+c\theta_0 y\left(\dfrac{1}{(x-r_{\epsilon}(t))^2+y^2}-\dfrac{1}{(x-\frac{R^2}{r_{\epsilon}(t)})^2+y^2} \right)\\
\theta_0x-c\theta_0\left(\dfrac{x-r_{\epsilon}(t)}{(x-r_{\epsilon}(t))^2+y^2}-\dfrac{x-\frac{R^2}{r_{\epsilon}(t)}}{(x-\frac{R^2}{r_{\epsilon}(t)})^2+y^2} \right)
\end{array}
\right).
\end{equation}
The previous vector is continuous for all $(t,(x,y);\epsilon)\in\R\times \R^2\times(0,\epsilon_0)$ whereas $(x,y)\notin \{(r_{\epsilon}(t),0),(\frac{R^2}{r_{\epsilon}(t)},0)\}$. Since $r_{\epsilon}(t)=r_0+o(\epsilon)$ and $\frac{R^2}{r_0}>R$ we can take $\epsilon$ small enough to ensure that $\frac{R^2}{r_{\epsilon}(t)}>R$. On the other hand, by Theorem~\ref{thm:retratos} $(b)-(d)$, $d_H(\mathcal A_M,r_0)>d_H(\mathcal A_M,\gamma)>0$ where $\gamma$ denotes the saddle connection that forms the outer boundary of the period annulus and $d_H$ denotes the Haussdorf distance of non-empty compact subsets of $\R^2$. Thus, by continuity, there exists $\epsilon_0>0$ small enough such that $d_H(\mathcal A_M,r_{\epsilon}(t))>d_H(\mathcal A_M,\gamma)>0$ for all $\abs{\epsilon}<\epsilon_0$. This implies that $\nabla\Psi(t,\eta;\epsilon)$ is continuous for all $(t,\eta;\epsilon)\in\R\times\mathcal A_M\times (0,\epsilon_0)$ as we desired. Moreover, for fixed $(t,\epsilon)\in \R\times(0,\epsilon_0)$, since $d_H(\mathcal A_M,r_{\epsilon}(t))>d_H(\mathcal A_M,\gamma)>0$ then $\nabla\Psi(t,\eta;\epsilon)\in C^1(\mathcal A_M)$. Consequently, $\nabla\Psi(t,\eta;\epsilon)$ is locally Lipschitz continuous in $\mathcal A_M$.
Then we can apply Theorem~\ref{thm:FSZ} and, particularly, Corollary~\ref{cor:FSZ} to show that for any positive integer $N$ there exists $0<\bar\epsilon_N<\epsilon_0$ such that if $\abs{\epsilon}<\bar\epsilon_N$ then system $\dot\eta=J\nabla\Psi(t,\eta;\epsilon)$ has at least $N$ periodic solutions in $\mathcal A_M\subset D_R$. Finally, by construction of $\mathcal A_M$ those periodic solutions have zero winding number with respect to the vortex.
\end{demo}

\section{Conclusions}

The main result of Section \ref{sec2} has a natural reading for the underlying physical model. Note that $\rho_0$ is the ratio between the path and domain radii respectively, while $\phi_0$ measures the relation between the path angular speed and the vortex strength. The sign of $\phi_0$ indicates if the sense of rotation of the vortex and the path is the same or opposite. For example, fixing the positive parameters $R,\Gamma,r_0$ and leaving $\theta_0$, for small positive $\theta_0$ there are no equilibria. Then, a first bifurcation point is $\theta^*_0=\frac{\Gamma}{2\pi R^2}\frac{R-r_0}{R+r_0}$, where a degenerate saddle appears at $(-R,0)$. A second bifurcation point appears at $\theta^{**}_0=\frac{\Gamma}{2\pi R^2}\frac{R+r_0}{R-r_0}$. For $\theta_0\in ]\theta^{*}_0,\theta^{**}_0[$, there is a center and two hyperbolic saddles in the border of the domain connected by an heteroclinic. They travel along the border until they collide at $\theta_0=\theta^{**}_0$ into a degenerate saddle, that enter into the domain as an hyperbolic saddle for values above $\theta^{**}_0$. On the other hand, for negative values of $\theta_0$, corresponding to opposite rotating sense of the vortex and the stirring protocol, we identify a typical saddle-node bifurcation. 

In the identified bifurcation scheme, saddles are connected by heteroclinic or homoclinic orbits that constitute barriers for the flux transport. Around the centers, the particles rotate with different periods, and this fact makes possible an application of a suitable result for perturbed hamiltonians, proving that for a perturbed vortex path there exist infinitely many periodic solutions that do not rotate around the vortex. The problem to identify more general classes of vortex protocols that generate this kind of periodic orbits with zero winding number is still open.

\end{document}